\documentclass[12 pt]{amsart}
\usepackage{amscd,amssymb,amsmath,amsthm, appendix}
\input xy
\xyoption{all}

\def \p {{\mathcal P}_{r,k}}
\def \gbar {{(g-1)}}
\newcommand{\comment}[1]{}

\newtheorem{theorem}{Theorem}
\newtheorem {lemma}{Lemma}

\newtheorem {proposition}{Proposition}

\theoremstyle{definition}
\newtheorem{remark}{Remark}

\theoremstyle {definition}

\begin{document}
\baselineskip=16pt

\title[$GL$ Verlinde numbers and the Grassmann TQFT]{$GL$ Verlinde numbers and the Grassmann TQFT}
\author {Alina Marian and Dragos Oprea}
\address {Department of Mathematics}
\address {University of Illinois at Chicago}
\email {alina@math.uic.edu}

\address {Department of Mathematics}
\address {University of California, San Diego}
\email{doprea@math.ucsd.edu}

\date{}
\maketitle

\section{Introduction}

These notes are concerned with moduli spaces of bundles on a smooth projective curve. Over them we consider determinant line bundles and their holomorphic Euler characteristics, the Verlinde numbers. The goal is to give a brief exposition of  the two-dimensional topological quantum field theory that captures the structure of the $GL$ Verlinde numbers, associated with spaces of bundles with varying determinant.  Our point of view is to emphasize the close connection with another TQFT, the quantum cohomology of the Grassmannian. 

Two different geometries are related here, the moduli of bundles on a curve $C$ and the space of maps from $C$ to a suitable Grassmannian. The connection between them was established in the classic paper \cite{witten} where the open and closed invariants of the $GL$ Verlinde TQFT, in all genera, were exhaustively written in both geometries. On the mathematical side, it was shown \cite{A} that the underlying algebras of the two TQFTs are isomorphic, as the genus zero three-point invariants match. Importantly however, the {\it metrics} of the associated Frobenius algebras are different: quantum cohomology uses the Poincar\'{e} pairing on the cohomology of the Grassmannian, while the Verlinde theory uses an intersection product on a space of higher degree maps from $C$ to the Grassmannian.  The TQFTs therefore turn up different invariants overall. The higher genus $GL$ Verlinde invariants, open or closed, have not been systematically written down in the mathematics literature although they were shown in \cite{witten} to have compelling closed-form geometric expressions. We found it useful therefore to render the results of \cite{witten} in standard mathematical language, also with a view toward future
studies of q-deformations of ordinary two-dimensional Yang Mills theory.

The exposition is organized as follows. The central point is presented in the final Section 5, where the Verlinde TQFT is explicitly written. Prior to this, we recall briefly the notion of a two-dimensional TQFT in the next section, then we introduce in our context, on a smooth projective curve $C$, the two spaces of interest:  the Grothendieck Quot scheme, and the moduli space of semistable bundles. We present the former here primarily as compactifying the space of maps from the curve to a Grassmannian. Relevant aspects of the geometry and intersection theory of the two spaces are discussed. The last section studies the relation between them, in the form of the $GL$ Verlinde TQFT, which we also refer to as the Grassmann TQFT.  
\vskip.1in
{\bf Acknowledgements.} The notes follow a series of lectures given by A. M. at the Geometry Summer School of the Instituto Superior T\'{e}cnico in Lisbon, in July 2009. She thanks the organizers Ana Cannas da Silva and Rui Fernandes for the warm hospitality and fantastic time in Lisbon during the School. The paper was written while A. M. was visiting Harvard in the fall of 2009. Partial support for both authors was provided by the NSF. 

\section{Generalities on two-dimensional TQFTs}

\label{tqftgen}

We consider the category ${\bf {2Cob}}$, in which 

\begin{itemize}
\item[(i)] the objects are one-dimensional compact oriented manifolds {\it i.e.}, finite unions of oriented circles; \vskip.03in
\item [(ii)] the morphisms are (diffeomorphism classes of) oriented cobordisms; \vskip.03in
\item [(iii)] composition of morphisms is concatenation of cobordisms; \vskip.03in
\item [(iv)] there is a tensor structure given by taking disjoint unions of objects. \vskip.03in
\end{itemize}

Let ${\mathbf {Vect}}_{\mathbb C}$ be the category of ${\mathbb C}$-vector spaces. A two-dimensional ${\mathbb C}$-valued TQFT is a symmetric monoidal functor 
$$F: \, \, {\bf {2Cob}} \, \longrightarrow \, {\mathbf {Vect}}_{\mathbb C}.$$  
There is a basic vector space $H$ in the theory, representing the value of the functor $F$ at the oriented circle $S^1$. In addition, $F$ associates to the empty manifold the vector space ${\mathbb C}$. 

The datum of the functor is equivalent to the structure of a commutative Frobenius algebra on $H$. By definition this comprises

\begin{itemize}
\item [(i)] a commutative associative multiplication $$H \otimes H \stackrel{{\bullet}} {\to} H$$ with identity element, and 

\item [(ii)] a symmetric nondegenerate pairing $$(\cdot  \, , \,\cdot ): H \otimes H \to {\mathbb C}$$ satisfying the Frobenius property $$(a\cdot b, \, c) = (a, \, b \cdot c).$$\end{itemize}
Indeed, if 
$W_s^t(g)$ is the genus $g$ cobordism with $s$ inputs and $t$ outputs, then 
\begin{enumerate}
\item[(i)]
$F(W_2^1(0)): H \otimes H {\to} H$ is the algebra multiplication, 
\item[(ii)]
$F(W_0^1(0)): {\mathbb C} \to H$ is the identity element, 
\item[(iii)] 
$F(W_2^0 (0))$ gives the pairing $(\cdot \, , \, \cdot ).$
\end{enumerate}

\vskip.1in

Viewed as a cobordism from the empty manifold to the empty manifold, a closed surface of genus $g$ corresponds under $F$ to a homomorphism from ${\mathbb C}$ to ${\mathbb C}$, thus to a number $F(g),$ $$ F(g) = F (W_0^0 (g)).$$ 

Let us assume that $H$ has a preferred basis, $$H = \bigoplus_{\lambda} {\mathbb C} e_{\lambda}.$$ The vector space $H^{\otimes s}$ has a basis $e_{\underline \lambda}$ indexed by multi-indices $\underline \lambda=(\lambda_1, \ldots, \lambda_s)$: $$e_{\underline \lambda}=e_{\lambda_1}\otimes \cdots \otimes e_{\lambda_s}.$$ We denote by $F(g)_{\underline{\lambda}}^{\underline{\mu}}\,$ the matrix entries of the cobordism homomorphism
$$F (W_s^t (g)): H^{\otimes s} \longrightarrow H^{\otimes t}$$ in this basis. We thus have
$$F(W_s^t (g)): \, \, e_{\underline{\lambda}} \mapsto  F(g)_{\underline{\lambda}}^{\underline{\mu}} \, \, e_{\underline{\mu}},$$ where ${\underline{\lambda}}$, $\underline{\mu}$ are multi-indices (with $s$ and $t$ components respectively). The TQFT is equivalent to the data of the numbers $F(g)_{\underline{\lambda}}^{\underline{\mu}}$ satisfying gluing rules which reflect the functoriality,
\begin{equation}
\label{tqft}
\sum_{\underline \mu} F(g_1)_{\underline{\lambda}}^{\underline{\mu}} \,\,  F(g_2)_{\underline{\mu}}^{\underline{\nu}} = F (g_1 + g_2 +  t - 1)_{\underline{\lambda}}^{\underline{\nu}}.
\end{equation}
Here $t$ is the number of components of the multi-index ${\underline{\mu}},$ which is summed over. 
\vskip.2in

\section{The Quot scheme $Q_C ({\mathbb G} (r,n), d)$}
\vskip.1in

Let $C$ be a smooth complex projective curve of genus $g$. We let $Q_C ({\mathbb G} (r,n), d)$ denote the Grothendieck Quot scheme parametrizing rank $n-r$ degree $d$ quotients of the rank $n$ trivial sheaf on $C$. A point in the Quot scheme is given by a short exact sequence
$$0 \rightarrow E \rightarrow {\mathcal O}_C\otimes \mathbb C^n \rightarrow F \rightarrow 0.$$ While the kernel sheaf $E$ is always locally free, the quotient $F$ is in general a sum
$$F = \overline{F} \oplus T,$$ with $\overline{F}$ locally free and $T$ a torsion sheaf supported at finitely many points of the curve $C$.

The quotients $F$ which {\it are} locally free form an open locus in $Q_C ({\mathbb G} (r,n), d)$, and can be regarded as degree $d$ maps $$f:C\to \mathbb G(r, n)$$ from $C$ to the Grassmannian ${\mathbb G} (r, n)$ of $r$ planes in ${\mathbb C}^n.$ The Quot scheme may be viewed as compactifying the space ${\text{Mor}}_d (C, {\mathbb G} (r, n))$ of degree $d$ maps to ${\mathbb G} (r,n)$: $${\text{Mor}}_d (C, {\mathbb G} (r, n))\hookrightarrow Q_C ({\mathbb G} (r,n), d).$$

\vskip.2in

\subsection{Examples} When $C = {\mathbb P}^1$ and $r=1$, the Quot scheme $Q_{{\mathbb P}^1} ({\mathbb P}^{n-1}, d)$ is the projectivized space of $n$ homogeneous degree $d$ polynomials in ${\mathbb C} [x,y],$  $$Q_{{\mathbb P}^1} ({\mathbb P}^{n-1}, d) \simeq {\mathbb P}^{n (d+1) - 1}.$$ 

In general, when $r=1$ and $C$ has arbitrary genus, 
$Q_{C} ({\mathbb P}^{n-1}, d)$ parametrizes exact sequences $$ 0 \rightarrow L \rightarrow {\mathcal O_C\otimes \mathbb C^n} \rightarrow Q \rightarrow 0$$ where $L$ is a line bundle of degree $-d$. Equivalently, dualizing such exact sequences, points in the space are degree $d$ line bundles $L^{\vee}$ on $C$ together with $n$ sections, not all zero: $$\mathcal O_C\otimes {\mathbb C^n}^{\vee}\to L^{\vee}.$$  Let $\text{Jac}^d (C)$ be the Picard variety of degree $d$ line bundles on $C$, and let $$\pi: {\text{Jac}}^d (C) \times C \to \text{Jac}^d (C)$$ be the projection. For $d$ sufficiently large, 
$d \geq 2g-1,$
the push forward $\pi_{\star} {\mathcal P}$ of the Poincar\'{e} line bundle $${\mathcal P} \rightarrow {{\text{Jac}}^d (C) } \times C$$ is locally free, and its fiber over $[L] \in {\text{Jac}}^d (C)$ is the space $H^0 (C, L)$ of sections of $L$.  In this case, 
$$Q_{C} ({\mathbb P}^{n-1}, d) \simeq {\mathbb P}( (\pi_{\star} {\mathcal P})^{\oplus n}) \to {\text{Jac}}^d (C).$$ 

\vskip.1in

Although for arbitrary $r$ the Quot scheme does not have such a simple description, it remains true that the space is well-behaved in the regime of large degrees $d$:

\begin{theorem}
\cite{bdw}
For $d >> r, n, g,$ the space $Q_C ({\mathbb G} (r, n), d)$ is irreducible, generically smooth, and has the expected dimension.
\end{theorem}

\vskip.2in 

\subsection{Structures}

As a fine moduli space, the Quot scheme carries a universal sequence
$$ 0 \rightarrow {\mathcal S} \to {\mathcal O}\otimes \mathbb C^n \to {\mathcal Q} \to 0 \, \, \, \text{on} \, \, \, Q_C ({\mathbb G} (r,n), d) \times C,$$ with the universal subsheaf ${\mathcal S}$ being locally free. The tangent sheaf to $Q_C ({\mathbb G} (r,n), d)$ is given as 
$${\mathcal T}Q_C ({\mathbb G} (r,n), d) \simeq {\text{Hom}}_{\pi} ({\mathcal S}, {\mathcal Q}),$$ where $$\pi: Q_C ({\mathbb G} (r,n), d) \times C \to Q_C ({\mathbb G} (r,n), d)$$ is the projection. The obstruction sheaf is ${\text{Ext}}^1_{\pi} ({\mathcal S}, {\mathcal Q}).$ The expected dimension is 
$$e = nd - r(n-r) (g-1)$$ by the Riemann-Roch formula.

The Chern classes of the universal subsheaf are natural to consider for the intersection theory of $Q_C ({\mathbb G} (r,n), d)$. Fixing a basis $$1, \delta_1, \ldots, \delta_{2g}, \omega$$ for the cohomology of the curve $C$, we write
$$ c_k ({\mathcal S}^{\vee}) = a_k \otimes 1 + \sum_{i=1}^{2g} b_k^i \otimes \delta_i + f_k \otimes \omega, \, \, \, \, 1 \leq k \leq r,$$
where $$a_k \in H^{2k} (Q_C ({\mathbb G} (r,n), d), {\mathbb C}), \, \, \, b_k^i \in H^{2k-1} (Q_C ({\mathbb G} (r,n), d), {\mathbb C}),$$  $$f_k \in H^{2k-2} (Q_C ({\mathbb G} (r,n), d), {\mathbb C}). $$
Note that 
\begin{equation}
f_k = \pi_{\star} c_k ({\mathcal S}^{\vee}),
\end{equation}
 while for $p \in C$ and $${\mathcal S}_p = \left. {\mathcal S} \right |_{Q_C ({\mathbb G} (r,n), d) \times \{p\}} ,$$ we have
\begin{equation}
\label{atop}
a_k = c_k ({\mathcal S}_p^{\vee}). 
\end{equation}
 
When $d$ is large so that $Q_C ({\mathbb G} (r, n), d)$ is irreducible, top intersections of the tautological $a, b$ and $f$ classes can be evaluated meaningfully against the fundamental class. For arbitrary degrees, the Quot scheme may be reducible and oversized. However, intersection theory can still be pursued in a virtual sense, by pairing Chern classes against a virtual fundamental cycle of the expected dimension, which the Quot scheme possesses:

\begin{theorem} \cite{cfk}, \cite{mo}
The Quot scheme $Q_C ({\mathbb G} (r,n), d)$ has a two-term perfect obstruction theory and a virtual fundamental class of expected dimension $$[Q_C ({\mathbb G} (r, n), d)]^{vir} \in A_e (Q_C ({\mathbb G} (r, n), d)).$$ 
\end{theorem}

\vskip.1in

{\it Proof.} We show that the tangent-obstruction complex for $Q_C ({\mathbb G} (r, n), d)$ admits a resolution
\begin{equation}
\label{tangentobs}
0 \to {\text{Hom}}_{\pi} ({\mathcal S}, {\mathcal Q} ) \to {\mathcal A}_0 \to {\mathcal A}_1 \to {\text{Ext}}^1_{\pi} ({\mathcal S}, {\mathcal Q}) \to 0,
\end{equation}
where the sheaves ${\mathcal A}_0$ and ${\mathcal A}_1$ are locally free. The virtual fundamental class is then standardly constructed as described in \cite{lt}, using the two vector bundles ${\mathcal A}_0, {\mathcal A}_1.$ 

The resolution is easily obtained as follows. Let ${\mathcal O} (1)$ be a degree $1$ line bundle on the curve $C$, and denote by ${\mathcal S} (m), {\mathcal Q} (m)$ the twists of the tautological sheaves by the pullback  of ${\mathcal O} (m)$ on $C$ to the product $Q_C ({\mathbb G} (r, n), d) \times C$.  Let $m$ be large enough so that 
$$R^1\pi_{\star} {\mathcal S} (m) = R^1 \pi_{\star} {\mathcal Q} (m) = 0,$$ 
and so that the evaluation map
$$\pi^{\star}\left (  R^0 \pi_{\star} {\mathcal S} (m) \right ) \to {\mathcal S} (m)$$ is surjective. The pushforward sheaves $R^0 \pi_{\star} {\mathcal S} (m), R^0 \pi_{\star} {\mathcal Q} (m)$ are then locally free. Further let ${\mathcal K}$ be the kernel
$$0 \rightarrow {\mathcal K} \to \pi^{\star} \left ( R^0 \pi_{\star} {\mathcal S} (m) \right ) \otimes {\mathcal O} (-m) \to {\mathcal S} \to 0.$$
Applying the functor $\text{Hom}_{\pi} (\cdot, {\mathcal Q} )$ gives
$$ 0 \to {\text{Hom}}_{\pi} ({\mathcal S}, {\mathcal Q}) \to \left ( R^0 \pi_{\star}{\mathcal S} (m) \right ) ^{\vee} \otimes R^0 \pi_{\star} {\mathcal Q} (m) \to {\text{Hom}}_{\pi} ({\mathcal K}, {\mathcal Q}) \to $$
$$\to {\text{Ext}}^1_{\pi} ({\mathcal S}, {\mathcal Q}) \to 0.$$  
Continuing this sequence one more term we get ${\text{Ext}}^1_{\pi} ({\mathcal K}, {\mathcal Q} ) = 0,$ so the sheaf $${\mathcal A}_1 =_{\text{def}} {\text{Hom}}_{\pi} ({\mathcal K}, {\mathcal Q})$$ is locally free. Also, $${\mathcal A}_0 =_{\text{def}} \left ( R^0 \pi_{\star}{\mathcal S} (m) \right )^{\vee} \otimes R^0 \pi_{\star} {\mathcal Q} (m),$$ is locally free.  
\qed

\subsection{Intersections} In this section, we will consider the (virtual) intersection theory of Quot schemes.

We start by pointing out the compatibility of the virtual fundamental class with the natural embedding, for $p \in C$,
\vskip.05in
$$\iota_p:  \, \, Q_C ({\mathbb G} (r, n), d) \, \, \hookrightarrow  \, \, Q_C ({\mathbb G} (r, n), d + r), $$ 
given by 
$$\{E \hookrightarrow {\mathcal O}_C\otimes \mathbb C^n  \} \, \, \, \mapsto \, \, \, \{ E(-p) \to E \to {\mathcal O}_C\otimes \mathbb C^n \}.$$ 
A degree $-d-r$ subsheaf $$E' \hookrightarrow {\mathcal O}_C\otimes \mathbb C^n$$ comes from $Q_C ({\mathbb G} (r, n), d)$ if the dual map $${{\mathcal O}_C\otimes \mathbb C^n}^{\vee}
\to {E'}^{\vee}$$ is zero at $p$. 
The image of the degree $d$ Quot scheme inside the degree $d+r$ space is therefore the zero locus of the dual universal map $${{\mathcal O}\otimes \mathbb C^n}^{\vee} \to {\mathcal S}_p^{\vee} \, \, \, \text{on} \, \, \, Q_C ({\mathbb G} (r, n), d +r). $$ This relationship is reflected on the level of the virtual fundamental classes for the two spaces. We recall that $a_r$ is the top Chern class of the universal subsheaf ${\mathcal S}_p^{\vee}$ before noting that 
\begin{proposition} \cite{mo} \label{consistent} The equality  
\begin{equation}\\
{\iota_p}_{\star} [Q_C ({\mathbb G} (r, n), d)]^{vir} = a_r^n \cap [Q_C ({\mathbb G} (r, n), d +r) ]^{vir}
\end{equation}\vskip.01in 
\noindent holds in $A_{\star} (Q_C ({\mathbb G} (r, n), d+r)).$
\end{proposition}

The intersection theory of $a$-classes is well understood. Top intersections are given in closed form by the {\it Vafa-Intriligator formula.} Furthermore, in the large-degree regime, the intersection numbers express counts of maps from the curve $C$ to the Grassmannian ${\mathbb G} (r, n),$ satisfying incidence constraints. More precisely, we have:

\begin{theorem}\label{bertramthm}
(i) \cite{bertram}, \cite{st}, \cite{mo} Let $J (x_1, \ldots, x_r)$ be the symmetric function  
$$J (x_1, \ldots, x_r) =  n^r \cdot {x_1^{-1} \cdots x_r^{-1}} \, {\prod_{1 \leq i < j \leq r} (x_i - x_j)^{-2}}.$$ Let $P (a_1, \ldots, a_r)$ be a top degree polynomial in the Chern classes of ${\mathcal S}_p^{\vee}.$ Then
\begin{equation}
\label{vaint}
\int_{[Q_C ({\mathbb G} (r, n), d)]^{vir}} P (a_1, \ldots , a_r) = u \cdot \sum_{\lambda_1, \ldots, \lambda_r} R (\lambda_1, \ldots , \lambda_r) \, J^{g-1} (\lambda_1, \ldots, \lambda_r),
\end{equation}
where $R$ is the symmetric polynomial obtained by expressing $P (a_1, \ldots a_r)$ in terms of the Chern roots of ${\mathcal S}_p^{\vee}$. The sum is taken over all $\binom{n}{r}$ tuples $$(\lambda_1, \ldots, \lambda_r)$$ of distinct $n$-roots of $1$. Here $$u=(-1)^{(g-1)\binom{r}{2}+d(r-1)}.$$
\text{(ii)} \cite{bertram} When $Q_C ({\mathbb G} (r, n), d)$ is irreducible of the expected dimension, the above intersection counts the number of degree $d$ maps from the curve $C$ to ${\mathbb G}(r, n)$ sending fixed distinct points of $C$ to special Schubert subvarieties of the Grassmannian, each Schubert variety matching an appearance of an $a$-class in the top monomial $P.$
\end{theorem}

The intersection numbers appearing in Theorem \ref{bertramthm} were written down in \cite{i}. Mathematical proofs have relied either on degenerations of the Quot scheme to genus zero, or on equivariant localization. Degeneration arguments use the enumerativeness of the $a$-intersections in the large-degree situation. 

By contrast, intersections involving $f$-classes do not give actual counts of maps, and explicit formulas for them have been relatively little explored. To describe one such formula, we let $$\sigma_i(x)=\sigma_i (x_1, \ldots, x_r) \text { and } \sigma_{i;k}(x)=\sigma_{i;k} (x_1, \ldots, x_r)$$ be the $i^{\text{th}}$ elementary symmetric functions in the variables $$x_1, \ldots, x_r \text{ and }x_1, \ldots,\widehat x_k, \ldots,  x_r$$ respectively. In the second set of variables, $x_k$ is omitted.

\begin{theorem} 
\cite{mo} Letting ${\mathcal D}_l, 2\leq l \leq r,$ be the first-order differential operator $${\mathcal D}_l = (g-1) (r-l+1)(n-r+l-1) \cdot \sigma_{l-1}(x)  + \sum_{k=1}^r \sigma_{l-1; k}(x) \,x_k \cdot \frac{\partial }{\partial x_k},$$ we have $$\int_{[Q_C ({\mathbb G} (r, n), d)]^{vir}} f_l \cdot P (a_1, \ldots, a_r) = \frac{u}{n} \sum_{\lambda_1, \ldots, \lambda_r} ({\mathcal D}_l R) (\lambda_1, \ldots, \lambda_r) \cdot J^{g-1} (\lambda_1, \ldots, \lambda_r).$$
\noindent The sum is over all $\binom{n}{r}$ tuples $(\lambda_1, \ldots, \lambda_r)$ of distinct $n$-roots of $1$.
\end{theorem}

It would be very interesting to generalize the Vafa-Intriligator formula to include all intersections of $f$ and $a$-classes.

We turn now to a discussion of the second geometry of interest. 

\section{The moduli space of semistable bundles}

\subsection{Basics} We consider vector bundles of rank $r$ and degree $d$ on the smooth curve $C$. We recall briefly the main facts in the moduli theory of semistable vector bundles on $C$. The family of {\it all} vector bundles of fixed topological type is not bounded, as one can immediately verify looking at vector bundles on ${\mathbb P}^1.$ A notion of stability is required to get a bounded problem.  For any vector bundle $E$, its {\it slope} $\mu (E)$ is defined as the ratio
$$\mu (E) = \frac{\text{degree} (E)}{\text{rank} (E)}.$$
A vector bundle $E$ is said {\it stable} ({\it semistable}) if for all subbundles $F \hookrightarrow E,$  $$\mu (F) < \mu (E) \, \, \, (\mu (F) \leq \mu (E)).$$
It follows easily that
\begin{lemma}
(i) If $E$ is semistable with $\mu (E) \geq 2g-1,$ then $H^1 (E) = 0.$ \\
(ii) If $E$ is semistable with $\mu (E) \geq 2g,$ then the evaluation map of sections $$H^0 (E) \otimes {\mathcal O_C} \to E$$ is surjective.
\end{lemma}

{\it Proof:} Indeed, by Serre duality, $H^1 (E)  \simeq H^0 (E^{\vee} \otimes K_C )^{\vee}, $ where $K_C$ denotes the canonical bundle. Let $L \hookrightarrow K_C$ be the image of an assumed nonzero homomorphism $\phi: E \to K_C$. $E$ is semistable and $L$ is a {\it quotient} of $E$, so we must have $$\mu (E) \leq \mu (L) = \deg (L) \leq \deg (K_C) = 2g-2.$$   
This contradicts the assumption, so there are no nonzero such homomorphisms and $H^1 (E) = 0.$  Regarding (ii), for any $p \in C,$ taking cohomology for the sequence  
$$ 0 \to E (-p) \to E \to E_p \to 0,$$ 
and using the vanishing of (i), it follows that the fiber of $E$ at $p$ is generated by global sections. 
\qed

\vskip.2in
Fixing a line bundle ${\mathcal O} (1)$ of degree $1$ on $C$, there is therefore an integer $m$ such that for all semistable rank $r$ and degree $d$ vector bundles $E$, we have 
$$H^1 (E(m)) = 0 \, \, \, \text{and} \, \, \, H^0 (E(m)) \otimes \mathcal O_C \to E (m) \to 0.$$ Any semistable $E$ can be thus realized as a quotient
$${\mathcal O}_C^{\oplus q} (-m) \to E \to 0, \, \, \, \text{with} \, \, \, q = \chi (E(m)),$$ 
{\it i.e.}, as a point in the Quot scheme $$\text{Quot}_{C}^{r,d} ({\mathcal O}_C^{\oplus q} (-m))$$ of quotients of ${\mathcal O}_C^{\oplus q} (-m)$ of rank $r$ and degree $d$. 
The group $SL(q)$ acts on this Quot scheme, with a standard linearization. On the locus of vector bundle quotients $E$ in $\text{Quot}_{C}^{r,d} ({\mathcal O}_C^{\oplus q} (-m))$ for which the quotient map induces an isomorphism 
 $$H^0 ({\mathcal O}_C^{\oplus q}) \simeq H^0 (E(m)), $$ stability in the geometric invariant theory sense coincides with slope stability.  
Restricting further to semistable quotients, we have an $SL(q)$-invariant open subscheme $${\text{Quot}}^{ss} \subset \text{Quot}_{C}^{r,d} ({\mathcal O}_C^{\oplus q} (-m)).$$ 
The GIT quotient

$${\text{Quot}}^{ss} // SL(q) =_{\text{def}} U_C (r, d)$$ 

\vskip.2in

\noindent{is} an irreducible normal projective variety of dimension $r^2 (g-1) + 1,$ the moduli space of semistable vector bundles of rank $r$ and degree $d$. The open subset $$U_C^s (r,d) \subset U_C (r,d)$$ parametrizing isomorphism classes of {\it stable} vector bundles is smooth and its complement has codimension at least 2 in $U_C (r,d).$  For details on this standard construction, we refer the reader to \cite{lepotier}.

\subsection{Line bundles on the moduli space and their Euler characteristics}

Twisting vector bundles by a line bundle of degree $1$ on $C$ gives an isomorphism $$U_C (r, d) \cong U_C (r, d+r),$$ so the dependence on degree is only modulo $r$. We assume further for simplicity that $$d=0.$$ All constructions can be easily duplicated in the arbitrary degree situation.

\vskip.1in

When $r = 1,$ we have $$U_C (1, 0) \simeq {\text{Jac}} (C),$$ the Picard variety of degree $0$ line bundles on $C$. Note that for a fixed line bundle $M$ on $C$ of degree $g-1$, 
$$\chi (L \otimes M) = 0 \, \, \text{for} \, \, L \in {\text{Jac}} (C).$$
The classical theta divisor relative to $M$ is defined as
$$\Theta_{1, M} = \{L \in {\text{Jac}} (C) \, \, \text{such that} \, \, h^0 (L \otimes M) \neq 0 \}.$$ 
Sections of the tensor powers of the line bundle ${\mathcal O} (\Theta_{1, M})$ are the classical theta functions, and 
\begin{equation}
\label{thetaclassic}
h^0 ( {\text{Jac}} (C), {\mathcal O} (k\, \Theta_{1, M})) = \chi ( {\text{Jac}} (C), {\mathcal O} (k\, \Theta_{1, M})) = k^g
\end{equation}
 is the dimension of the space of level $k$ theta functions.

 For $r > 1,$ we have similarly, when  $M$ is as before a line bundle of degree $g-1$ on $C$, 
$$\chi (E \otimes M) = 0 \, \, \text{for} \, \, E \in U_C (r, 0),$$ 
and we set

 \begin{equation}
\label{thetadiv}
 \Theta_{r, M} = \{ E \in U_C (r, 0) \, \, \text{such that} \, \, h^0 (E \otimes M) \neq 0 \}.
 \end{equation}

\vskip.2in

\noindent{As} in the $r = 1$ case in fact, the divisor $\Theta_{r, M}$ has a determinantal  scheme structure: for a family $${\mathcal E} \to S \times C$$ of semistable rank $r$ degree 0 vector bundles, flat over $S$, we consider a resolution 
$$ 0 \to R^0\pi_{\star}( {\mathcal E} \otimes p_C^{\star} M) \to {\mathcal F}_0 \stackrel{\varphi}{\to}  {\mathcal F}_1 \to R^1 \pi_{\star} ({\mathcal E} \otimes p_C^{\star} M )\to 0 $$
of the direct image complex $$R\pi_{\star}( {\mathcal E} \otimes p_C^{\star}M),$$  so that ${\mathcal F}_0, {\mathcal F}_1$ are locally free. Here we denoted by $\pi$ and $p_C$ the projections 
$$S \times C \stackrel{\pi}{\to} S, \, \, \, S \times C \stackrel{p_C}{\to} C.$$ The pullback of $\Theta_{r, M}$ to $S$ is then the degeneracy locus of $\varphi$. 
The line bundle ${\mathcal O} (\Theta_{r, M})$ is the descent of the determinant line bundle $$\det R\pi_{\star} ({\mathcal E} \otimes p_C^{\star} M)^{-1}$$ from the Quot scheme $\text{Quot}_{C}^{r,d} ({\mathcal O}_C^{\oplus q} (-m))$, with ${\mathcal E}$ being the universal quotient.

\vskip.2in

The Picard group of $U_C (r, 0)$, described in \cite{dn}, is generated by the theta line bundles ${\mathcal O} (\Theta_{r, M})$ as $M$ varies in $\text{Pic}^{g-1} (C).$ 
To state this precisely, let $$\det: U_C (r, 0) \to \text{Jac} (C)$$ be the morphism sending bundles to their determinants. The following holds.
\begin{theorem} \cite{dn} 
\comment{
Identify also, as is standard,
\begin{equation}
\label{picc}
{\text{Pic}}^0 (\text{Jac} (C)) \cong \text{Jac} (C).
\end{equation}}
(i)  Consider $$\iota: SU_C (r, \mathcal O) \hookrightarrow U_C (r, 0)$$ the moduli space of bundles with trivial determinant. 
The restriction $${\mathcal L} = _{\text{def}} \iota^{\star} {\mathcal O} (\Theta_{r, M}),  \, \, \, \, \, \, $$ is independent of the choice of $M$ in ${\text{Pic}}^{g-1} (C)$ and 
$${\text{Pic}} (SU_C (r, \mathcal O)) \cong {\mathbb Z} \, {\mathcal L}.$$
(ii) $${\text{Pic}}\,(U_C (r, 0)) \cong {\mathbb Z} \, {\mathcal O} (\Theta_{r, M}) \oplus {{\det}^{\star}} (\text{Pic}^0 ( {\text{Jac}} (C))).$$ \comment{where the identification \eqref{picc} was made.  }
\end{theorem}

As in the classical case, the theta bundles on $U_C (r, 0)$ and $SU_C (r, \mathcal O)$ have no higher cohomology, so their holomorphic Euler characteristics give also the dimension of their spaces of sections. Explicit expressions for them, known as {\it Verlinde formulas}, were derived by several methods, and are significantly more complicated than \eqref{thetaclassic}. The formulas are very similar for $k$ powers of ${\mathcal L}$ on $SU_C (r, \mathcal O)$ and of ${\mathcal O} (\Theta_{r, M})$ on $U_C (r, 0). $ A slightly simpler and more convenient  expression arises however for the twist 
$${\mathcal O} (k \, \Theta_r) \otimes  {\det}^{\star} {\mathcal O} (\Theta_1) \in {\text{Pic}} (U_C (r, 0)).$$ \vskip.03in
\noindent{Here we suppressed reference degree $g-1$ line bundles for the theta bundles, as the holomorphic Euler characteristic is independent of these choices.} Writing also, to simplify notation,  $\Theta_r$ and $\Theta_1$ for the line bundles ${\mathcal O} (\Theta_r)$ and ${\mathcal O} (\Theta_1)$, we have  
\begin{eqnarray}
\label{verlinde}
V_g^{r, k} & =_{\text{def}} & h^0 (U_C (r, 0), \Theta_r^k  \otimes {\det}^{\star} \Theta_1 ) =  \chi (U_C (r, 0), \Theta_r^k \otimes {\det}^{\star} \Theta_1)\\  & = & \sum_{\stackrel{S \sqcup T = \{1,  \ldots, r+k \}}{ |S| = r}} \prod_{\stackrel{s\in S}{ t\in T}} \left | 2 \sin \pi \frac{s-t}{r+k} \right |^{g-1}. \nonumber
\end{eqnarray}
The sum is over the $\binom{r+k}{r}$ partitions of the first $r+k$ natural numbers into two disjoint subsets $S$ and $T$ of cardinalities $r$ and $k$.  Note that the numbers $V_g^{r, k}$ depend solely on the genus $g$ of $C$, the rank $r$, and the level $k$.

\vskip.3in

\subsection{Parabolic counterparts}

We would like to formulate degeneration rules for the Verlinde numbers $V_g^{r,k}.$ To this end, we turn to decorated moduli spaces of rank $r$ vector bundles on $C$. In addition to $r$, we think of the {\it level} $k$ as fixed.  

We denote by $\p$ the set of Young diagrams with at most $r$ rows and at most $k$ columns. Enumerating the lengths of the rows, we write a diagram $\lambda$ as  $$\lambda = (\lambda^{1}, \ldots, \lambda^{r}), \, \, k \geq \lambda^{1} \geq \cdots \geq \lambda^{r} \geq 0.$$ Such vectors can also be regarded as highest weights for irreducible representations of the unitary group $U(r)$, bounded by $k$. Further, the lengths of {\it columns} in a partition $\lambda \in \p$ give a flag type on an $r$-dimensional vector space. Let $\text{Fl}_{\lambda}$ be the corresponding flag variety, with associated Borel-Weil line bundle ${\mathcal N}_{\lambda}.$ 

We  consider the curve $C$ together with a finite set $I$ of distinct points on it, and partitions $\lambda_p \in \p$ labeled by the points $p \in I.$ 
A vector bundle $E$ together with a choice of a flag in each of its fibers over the points in $I$, $$0\subset E_{1,p}\subset E_{2, p}\subset \ldots \subset E_{k, p} \subset E_p$$  with flag type given for each $p \in I$ by the partition $\lambda_p,$ 
is referred to as a {\it parabolic} vector bundle of type  $\underline{\lambda} = (\lambda_p)_{p \in I}$. The lengths of {\it rows} in a partition $\lambda_p$ add the datum of a set of weights to the flag type at $p$, and define a parabolic slope for $E$,  
\begin{equation}
\label{parslope}
\mu_{\text{par}} (E) = \frac{d}{r} + \frac{|\underline\lambda|}{rk},
\end{equation}
 with $|\underline \lambda|$ being the total number of boxes in all partitions $\lambda_p, p \in I.$ As in the case of undecorated bundles, the slope comes with a notion of semistability, and there is a coarse projective moduli space $U_C (r, d, \underline{\lambda})$ of semistable rank $r$ degree $d$ parabolic vector bundles 
of type ${\underline{\lambda}}$, introduced in \cite{ms}. 

The construction is similar to that of the undecorated space $U_C (r, d).$ Its brief description here follows \cite{pauly}. To start, let $\Omega$ be the open locus in the Quot scheme $\text{Quot}_{C}^{r,d} ({\mathcal O}_C^{\oplus q} (-m))$  where the universal quotient sheaf $${\mathcal Q}\,  \to \, \text{Quot}_{C}^{r,d} ({\mathcal O}_C^{\oplus q} (-m)) \, \times \, C$$ is locally free, and in addition each quotient $$\mathcal O_C^{\oplus q} (-m) \to E$$ in $\Omega$ gives an isomorphism $$H^0 ({\mathcal O}_C^{\oplus q}) \simeq H^0 (E(m)).$$ 
For each point $p \in I$, consider next the restriction $${\mathcal Q}_p = \left . {\mathcal Q} \right |_{\Omega \times \{p\}}$$ of the universal quotient bundle, and its associated flag bundle ${Fl}_{\lambda_p},$ where the flag type is specified by the partition $\lambda_p$. Let $R$ be the product over $\Omega$ of the flag bundles for each $p \in I$,
$$ R = {{Fl}}_{\lambda_{p_1}} \times_{\Omega} \cdots \times_{\Omega} {{Fl}}_{\lambda_{p_n}}.$$
The moduli space of semistable parabolic vector bundles of type $\underline{\lambda}$ is the GIT quotient 

$$U_C (r, d, \underline{\lambda})  =_{\text{def}} R^{ss}//SL(q),$$ 

\vskip.2in

\noindent{where $R^{ss}$ is the open semistable locus in $R$ defined in terms of the slope \eqref{parslope}.}

We describe natural theta bundles over $U_C(r, d, \underline \lambda)$. One can consider on $\Omega$ the level $k$ determinant line bundle 

$$\left ( \det  R\pi_{\star}({\mathcal Q})\right ) ^{- k},$$ 
where as usual $$\pi: \Omega \times C \to \Omega$$ is the projection. Furthermore each flag bundle $Fl_{\lambda_{p}}$ carries a natural line bundle $${\mathcal N}_p\to Fl_{\lambda_p}$$ restricting  fiberwise over $\Omega$ to the Borel-Weil ample line bundle ${\mathcal N}_{\lambda_p}$ on the flag variety ${\text{Fl}}_{\lambda_p}.$ More precisely, the line bundle ${\mathcal N}_p$ is the tensor product of the determinants of all universal quotients on the flag bundle $Fl_{\lambda_p}.$   
Under the condition
\begin{equation}
\label{selectionrule}
kd+|\underline{\lambda} | \equiv 0 \mod r,
\end{equation}
the tensor product 
\begin{equation}\label{thetapar}\left ( {\det R\pi_{\star} ({\mathcal Q})}\right ) ^{-k} \bigotimes_{p\in I}{\mathcal N}_{p}\otimes (\det Q_{x})^{e}\end{equation} descends to a line bundle $${\mathcal L}_{\underline{\lambda}} \, \, \to \, \, U_C (r, d, \underline{\lambda})$$ on the GIT quotient. Here $x$ is a point on the curve (which will be omitted from the notation), and $$e=\frac{kd+|\underline{\lambda}|}{r}+k(1-g).$$ 
When $\underline{\lambda}$ consists of empty partitions, and $d=0$, we recover the space $U_C (r, 0)$ and the line bundle $\Theta_{r, M}^k$ where $M=\mathcal O((g-1)x).$
\vskip.1in

We set
\begin{equation}
V_{g,d}^{r, k} (\underline{\lambda}) = h^0 (U_C (r, d, \underline{\lambda}), {\mathcal L}_{\underline{\lambda}} \otimes {\det}^{\star} \Theta_1) = \chi (U_C (r, d, \underline{\lambda}), {\mathcal L}_{\underline{\lambda}} \otimes {\det}^{\star} \Theta_1).
\end{equation}\vskip.03in
\noindent The parabolic Verlinde numbers $V_{g, d}^{r,k} (\underline{\lambda})$ are given by explicit elementary formulas similar to \eqref{verlinde}. As they are generally elusive in the literature, we write them down in detail in the appendix, in the course of describing the relationship between $V_{g, d}^{r,k}(\underline \lambda)$ and intersections on the Quot scheme.

\section{The $GL$ Verlinde TQFT at fixed rank and level}

\subsection {Euler characteristics and intersections on the Quot scheme}

The theory of Euler characteristics of determinant line bundles over the moduli space $U_C (r, 0)$ is naturally related to the intersection theory of the space $${\text{Mor}}_d (C, {\mathbb G} (r, k+r))$$ of degree $d$ maps to ${\mathbb G} (r, k+r),$ 
where $$d \equiv 0 \mod r.$$
One of the most concrete aspects of this connection is the following remarkable formula \cite{witten} for the undecorated Verlinde numbers. Recall the top Chern class $a_r$, defined in \eqref{atop}, on the Quot scheme $Q_C ({\mathbb G} (r, k+r), d)$ compactifying  ${\text{Mor}}_d (C, {\mathbb G} (r, k+r))$. We define the integer $$t = \frac{d}{r} (k+r) - k (g-1),$$ so that the expected dimension of $Q_C(\mathbb G(r, k+r), d)$ equals $rt$. The Verlinde number $V_g^{r, k}$ can be expressed as  a top intersection 

\begin{equation}
 \label{remarkable} 
 V_g^{r, k} = \int_{[Q_C ({\mathbb G} (r, k+r), d)]^{vir}} a_r^t. 
\end{equation} 

\vskip.2in

\noindent Note that although $d$ is arbitrary divisible by $r$, Proposition \ref{consistent} ensures that \eqref{remarkable} gives the same answer for different values of $d$. 

\vskip.1in

It can be easily checked in fact that \eqref{remarkable} holds: the Vafa-Intriligator sum giving the right-hand side integral can be immediately written as the elementary formula \eqref{verlinde}. 
More satisfyingly, geometric arguments \cite{mo2} relate the intersection theory of  the space $U_C (r, d)$ with that of the Quot scheme $Q_C ({\mathbb G} (r, n), d)$ in the large $n$ limit. The particular expression of the Todd class appearing in holomorphic Euler characteristic calculations then recasts the Verlinde number $V_g^{r,k}$ as the intersection \eqref{remarkable} on the finite Quot scheme  $Q_C( {\mathbb G} (r, k+r), d).$ 

More generally, the decorated degree $d$ Verlinde number $$V_{g,d}^{r, k} (\underline{\lambda}) =\chi (U(r, d, \underline \lambda), \mathcal L_{\underline\lambda}\otimes {\det}^{\star}\Theta_1)$$ can be expressed as an intersection on the space of maps from $C$ to $\mathbb G(r, k+r)$ subject to incidence conditions with suitable Schubert cycles determined by the multipartition $\underline \lambda$. 

To explain the result, we need more notation. To an individual partition $\lambda \in {\mathcal P}_{r, k}$ we associate the Schur polynomial in the Chern roots $x_1, \ldots, x_r$ of the rank $r$ universal sheaf $\mathcal S_p^{\vee}$: $$\sigma_{\lambda} (x_1, \ldots, x_r) = \frac{\det (x_i^{\lambda^j + r-j})} {V (x_1, \ldots , x_r)},$$ where $V(x_1, \ldots, x_r)$ is the Vandermonde determinant. 
We denote the ensuing class $$a_{\lambda} = \sigma_{\lambda} ({\mathcal S}_p^{\vee}).$$ For a multipartition $\underline{\lambda} = (\lambda_1, \ldots, \lambda_n),$ we set $$a_{\underline \lambda}=a_{\lambda_1}\cdots a_{\lambda_n}.$$ Next, to a partition $$\lambda: k \geq \lambda_1 \geq \cdots \geq \lambda_r \geq 0 \, \, \text{in} \, \,  {\mathcal P}_{r, k},$$ we associate the conjugate partition $\lambda^{\star} \in {\mathcal P}_{r, k},$ $$ \lambda^{\star}:  k \geq k-\lambda_r \geq \cdots \geq k - \lambda_1 \geq 0.$$ The definition extends naturally to multipartitions $\underline \lambda.$

Under the assumption \begin{equation}\label{blah}kd+|\underline \lambda|\equiv 0\mod r,\end{equation} we have the equality
\begin{equation}
\label{rgeneral}
V_{g,d}^{r, k} (\underline{\lambda}) = \int_{[Q_C ({\mathbb G} (r, k+r), d)]^{vir}} a_{\underline \lambda^{\star}} \cdot a_r^t. 
\end{equation}
\vskip.05in
\noindent Here $t$ is taken to satisfy the dimension equation
$$|\underline\lambda^{\star}| + r t = (k+r) d - rk (g-1),$$ which is always possible when \eqref{blah} is satisfied. 

In degree $0$, formulas related to \eqref{rgeneral} were explicitly written down in \cite{oudompheng} in the process of establishing a level-rank duality on moduli of parabolic bundles.  The parabolic Verlinde numbers for arbitrary degree $d$ have been less explored, hence we will give an argument establishing \eqref{rgeneral} in the appendix.

\subsection{The Grassmann TQFT}

The Verlinde numbers $V_g^{r, k}$ are the closed invariants  
$$F(g) =  V_g^{r, k}$$ in a TQFT which we now describe. We refer to this theory equally as the $GL$ Verlinde, or the Grassmann TQFT. The theory was introduced in \cite{witten}, which we follow closely, while expressing the main facts in standard mathematical form. The fundamental vector space of the TQFT, together with a preferred basis, is   
$$H = \bigoplus_{\lambda \in \p} {\mathbb C} \lambda.$$ Considering the Grassmannian $G(r, k+r)$ and its tautological sequence
$$0 \rightarrow S \rightarrow {\mathcal O}\otimes \mathbb C^{r+k} \rightarrow Q \rightarrow 0,$$ we think of $x_1, \ldots, x_r$ as being the Chern roots of the dual tautological bundle $S^{\vee}$. In this case, the Schur polynomials $\sigma_{\lambda} (x_1, \ldots, x_r)$ give a basis for the cohomology of the Grassmannian, and we may view
$$H = \bigoplus_{\lambda \in \p} {\mathbb C} \, \sigma_{\lambda} = H^{\star} (G(r, k+r), {\mathbb C}).$$ 

The numbers $F(g)$ were written in the previous subsection as intersections on a suitable Quot scheme. The general matrix elements of $F(W_s^u (g))$ are integrals on the Quot scheme as well. We consider the Quot schemes for all degrees at once, setting $$Q_{C, r, k}=\coprod_{d} Q_C(\mathbb G(r, k+r), d).$$ As explained in the previous subsection, they come equipped with natural cohomology classes $a_{\underline \lambda}$, indexed by multipartitions. To start, for $\underline{\lambda}$ a multipartition with $s$ components, we define the matrix elements $F(g)_{\underline{\lambda}}$ of the homomorphism $$F(W_s^0 (g)): H^{\otimes s} \to {\mathbb C}$$ by
\begin{equation}
\label{partialmatrix}
F(g)_{\underline \lambda}=\int_{[Q_{C, r,k}]^{vir}} a_{\underline \lambda}\cdot a_r^{rg+k}.
\end{equation} We define the matrix elements $F(g)_{\underline \lambda}^{\underline \mu}$ in full generality by 

\begin{equation}\label{sum}F(g)_{\underline \lambda}^{\underline \mu}=\int_{[Q_{C, r, k}]^{vir}} a_{\underline \lambda}\cdot a_{\underline \mu^{\star}} \cdot a_r^{r(g+u)+k},
\end{equation} 

\vskip.18in

\noindent where $u$ is the number of components of the multipartition $\underline \mu$. Note that only one summand contributes to the infinite sum \eqref{sum}, since integration occurs only over the Quot scheme of degree \begin{equation} \label{d}d=\frac{|\underline \lambda|-|\underline \mu|}{k+r}+r(g+u).\end{equation} If this expression does not yield an integer {\it i.e.}, \begin{equation}\label{divisibil}|\underline \lambda|\not \equiv |\underline \mu| \mod k + r \end {equation} the matrix element $F(g)_{\underline \lambda}^{\underline \mu}$ is $0$. Letting $\mu$ in \eqref{sum} consist of no partitions, we recover \eqref{partialmatrix}. 

In the last subsection we show that the numbers $F(g)_{\underline{\lambda}}^{\underline{\mu}}$ satisfy the requisite gluing formula \eqref{tqft} of a TQFT. 

\subsection{Formulation in terms of Verlinde data.} The closed invariants $F(g)$ coincide with the undecorated Verlinde numbers $V_{g}^{r,k}$. Indeed, when $\underline\lambda$ and $\underline \mu$ both consist of no partitions, we obtain $$F(g)=\int_{[Q_{C, r, k}]^{vir}} a_r^{rg+k} =V_{g}^{r,k},$$ which is a  particular case of equation \eqref{remarkable} for $d=rg$. 

More generally, all matrix elements $F(g)_{\underline \lambda}^{\underline \mu}$ can be expressed as Verlinde numbers. Assuming $$|\underline \lambda|\equiv |\underline \mu|\mod r+k,$$ the degree given by \eqref{d} $$d=\frac{|\underline \lambda|-|\underline \mu|}{r+k}+r(g+u)$$ satisfies \eqref{blah}. Hence, by \eqref{rgeneral}, we have 
\begin{equation}\label{fglm}
F(g)_{\underline \lambda}^{\underline \mu}=\chi (U_C(r, d, \underline \lambda^{\star}, \underline \mu), \mathcal L_{\underline \lambda^{\star}, \underline {\mu}}\otimes {\det}^{\star}\Theta_1),
\end{equation}\vskip.03in
\noindent In particular, \begin{equation}\label{fgl}
F(g)_{\underline \lambda}=\chi \left(U_C(r, d, \underline \lambda^{\star}), \mathcal L_{\underline \lambda^{\star}} \otimes {\det}^{\star} \Theta_1\right).
\end{equation}

\begin{remark} {\it Comparison with the quantum cohomology of ${\mathbb G} (r, k+r).$}
There is a slight asymmetry between the roles of $\underline{\lambda}$ and $\underline{\mu}$ in \eqref{sum}, with only the number of components of the multi-index $\underline{\mu}$ appearing explicitly in the defining integral. This reflects a twist in the metric $F(W_2^0 (0))$ on the Frobenius algebra $H$. The metric is given by $$(\sigma_{\lambda}, \sigma_{\mu})=F(0)_{\lambda, \mu}=\int_{[Q_{{\mathbb P}^1, r,k}]} a_{\lambda} \,a_{\mu}\cdot a_r^{k},$$ which manifestly differs  
from the usual Poincar{\'{e}} pairing $$\int_{{\mathbb G} (r, k+r)}  a_{\lambda}\, a_{\mu}.$$ Turning now to the algebra structure on $H$, we have
$$\sigma_{\lambda}\cdot \sigma_{\mu}=\sum_{\nu} F(0)_{\lambda, \mu}^{\nu} \sigma_{\nu},$$ where $$F(0)_{\lambda, \mu}^{\nu}= \int_{[Q_{{\mathbb P^1}, r, k}]} a_{\lambda} \, a_{\mu} \, a_{\nu^{\star}}\cdot a_r^{k+r}= \int_{[Q_{{\mathbb P^1}, r, k}]} a_{\lambda} \, a_{\mu} \, a_{\nu^{\star}}.$$
The last integral gives precisely the structure constants of the quantum multiplication on $H^{\star} ({\mathbb G} (r, k+r), {\mathbb C})$ in the Schur basis. Therefore, we obtain an {\it algebra} isomorphism with quantum cohomology $$H \cong QH^{\star}({\mathbb G} (r, k+r)).$$ 
Being based on the Poincar{\'{e}} metric, the quantum cohomology as a TQFT is different however from the Grassmann TQFT given by the numbers $F(g)_{\underline{\lambda}}^{\underline{\mu}}$. This is accounted for by the disparity between the two metrics. 
\end {remark}
\vskip.1in

\begin{remark} {\it Comparison with the $SU(r)$ level $k$ fusion algebra.} A closely related theory is the well-studied $SL$ Verlinde TQFT described in \cite {beauville0} \cite {tuy}. The underlying vector space $$\widetilde H=\bigoplus_{\rho} \mathbb C\rho$$ is labeled by heighest weight representations $\rho$ of $SU(r)$ at level $k$. Most concretely, we think of $\rho$ as equivalence classes of partitions $\lambda\in \mathcal P_{r, k},$ where $$\lambda \sim \mu$$ if $\lambda$ and $\mu$ are obtained from one another by adding or subtracting the same number of boxes from all rows. 

In this basis, the matrix elements $\widetilde F(g)_{\underline \lambda}^{\underline \mu}$ of the theory are given as Verlinde numbers
$$\widetilde F(g)_{\underline \lambda}^{\underline \mu}=\chi(\mathcal L_{\underline \lambda, \underline \mu^{\star}})$$ where $$\mathcal L_{\underline \lambda, \underline \mu^{\star}}\to SU_C(r, \underline \lambda, \underline \mu^{\star})$$ is the level $k$ determinant bundle over the moduli space of parabolic bundles with trivial determinant. The degeneration formulas, known as factorization rules, were proved in \cite {tuy} using the connection with conformal blocks. 

The underlying algebra of the theory $\widetilde F$ admits a presentation which is similar, but not identical, to that of the quantum cohomology of $\mathbb G(r, k+ r)$ obtained in \cite {st}. The precise result is explained in Theorem $1.3$ of \cite {S}. 

\end{remark}

\begin {remark} {\it Level-rank duality.} Interchanging the rank $r$ and level $k$, we obtain a second TQFT from the enumerative geometry of the Grassmannian $\mathbb G(k, r+k)$. 
The vector space $\widehat H$ of the theory has a basis indexed by $k\times r$ partitions $$\widehat H=\bigoplus_{\lambda\in \mathcal P_{k,r}} \mathbb C\lambda.$$ The natural map $$H\to \widehat H,\,\,\, \lambda \to \lambda^{t}$$ sending a partition to its transpose is an isomorphism of TQFTs. Indeed, the equality of the matrix elements \begin{equation}\label{equality}F(g)_{\underline \lambda}^{\underline\mu}=\widehat F(g)_{\underline \lambda^{t}}^{\underline \mu^{t}}\end{equation} follows immediately from the explicit formula \eqref{dimm} of the appendix. In particular, the closed invariants satisfy the symmetry 
$$\chi\left(U_C(r, 0), \Theta_{r, M}^{k}\otimes {\det}^{\star} \Theta_{1, M}\right)=\chi\left(U_C(k, 0), \Theta_{k, M}^{r} \otimes {\det}^{\star} \Theta_{1, M}\right).$$ \vskip.03in
\noindent This equality reflects a geometrically induced isomorphism
$$H^0\left(U_C(r, 0), \Theta_{r, M}^{k}\otimes {\det}^{\star} \Theta_{1, M}\right)^{\vee}\cong H^{0}\left(U_C(k, 0), \Theta_{k, M}^{r} \otimes {\det}^{\star} (-1)^{\star} \Theta_{1, M}\right),$$ \vskip.05in \noindent 
proved in \cite {sd}. The symmetry of the open invariants \eqref{equality} on parabolic moduli spaces was similarly explained in \cite {oudompheng}. 

\end {remark}
\subsection {Degeneration rules} To prove that the matrix elements $F(g)_{\underline{\lambda}}^{\underline{\mu}}$ satisfy \eqref{tqft}, we show the two degeneration formulas
\begin{equation}
\label{deg1}
F(g)_{\underline \lambda}^{\underline \mu}=\sum_{\rho \in {\mathcal P}_{r, k}} F(g-1)_{\underline \lambda, \, \rho}^{\underline \mu, \, \rho},
\end{equation}
and 
\begin{equation}
\label{deg2}
F(g) _{\underline{\lambda}}^{\underline \mu} = \sum_{\rho\in {\mathcal P}_{r, k}} F(g_1)_{\underline{\lambda_1}, \,\rho}^{\underline \mu_1} \cdot F(g_2)^{\underline{\mu_2}, \, \rho}_{
\underline \lambda_2}
\end{equation}
for splittings $$g = g_1 + g_2,\,\, \underline \lambda=\underline \lambda_1+ \underline \lambda_2, \,\,\, \underline \mu=\underline \mu_1+\underline \mu_2.$$

The argument is standard. Suppose first that a smooth curve $C$ of genus $g$ degenerates to a nodal irreducible curve $C_0$ with one node $s$, and let $\widetilde{C}$ be the smooth genus $g-1$ curve normalizing $C_0$. We write the class of the diagonal $$\Delta \, \subset \,  {\mathbb G} (r, k+r) \times {\mathbb G} (r, k+r)$$ as 
$$[\Delta ] = \sum_{\rho \in {\mathcal P}_{r, k} } \sigma_{\rho} (x_1, \ldots, x_r) \, \sigma_{\rho^{\star}} (x_1', \ldots, x_r'), $$ where the primed variables are the Chern roots of the tautological bundle $S^{\vee}$ on the second Grassmannian. 
For any top polynomial $P (a_1, \ldots, a_r)$ and sufficiently large degrees $d$, the arguments of \cite{bertram} show that 
\begin{equation}
\label{deg3}
 \int_{Q_C ({\mathbb G} (r, k+r), d)} P(a_1, \ldots, a_r) = \sum_{\rho \in {\mathcal P}_{r, k}} \int_{Q_{\widetilde{C}} ({\mathbb G} (r, k+r), d)} P (a_1, \ldots, a_r)\, a_{\rho} \,a_{\rho^{\star}}. 
 \end{equation}
The equation expresses the fact that the space of maps ${\text{Mor}}_d (C_0, {\mathbb G} (r, k+r))$ is embedded in the larger space ${\text{Mor}}_d ({\widetilde{C}}, {\mathbb G} (r, k+r) )$ as the locus of maps $f$ with $f(s_1) = f(s_2),$ where the two points $s_1$ and $s_2$ lie over the node of $C_0.$ Thus
$${\text{Mor}}_d (C_0, {\mathbb G} (r, k+r)) = ev_2^{-1} (\Delta),$$
with $ev_2$ denoting the evaluation map $$ev_2: {\text{Mor}}_d ({\widetilde{C}}, {\mathbb G} (r, k+r)) \, \to {\mathbb G} (r, k+r) \times {\mathbb G} (r, k+r) $$ at $s_1$ and $s_2$. Intersections of $a$-classes on Quot compactifications of spaces of maps are moreover enumerative in the large degree regime, yielding \eqref{deg3}. Proposition \ref{consistent} then ensures that \eqref{deg3} holds in arbitrary degree when the integrals are evaluated against the virtual fundamental class.  

If we let $C$ degenerate to a reducible nodal curve with one node and two smooth irreducible components $C_1$ and $C_2$ of genera $g_1$ and $g_2$, such that $$g = g_1 + g_ 2,$$ a similar argument shows  
\begin{eqnarray}
\label{deg4}
\int_{[Q_C ({\mathbb G} (r, k+r), d)]^{vir}} P\cdot  Q (a_1, \ldots a_r) &=& \sum_{\rho \in {\mathcal P}_{r, k}} \sum_{d_1 + d_2 = d} \int_{[Q_{C_1} ({\mathbb G} (r, k+r), d_1)]^{vir} } P(a_1, \ldots, a_r)\, a_{\rho}  \nonumber \\ &\cdot&\int_{[Q_{C_2} ({\mathbb G} (r, k+r), d_2)]^{vir}} Q(a_1, \ldots, a_r)\,  a_{\rho^{\star}}.
\end{eqnarray}
Equation \eqref{deg4} is also argued geometrically in the large degree regime, where the intersections involved are enumerative. The passage to arbitrary degree and the virtual fundamental class is again via Proposition \ref{consistent}. 

The degeneration rule \eqref{deg1} follows from \eqref{deg3} taking $$P (a_1, \ldots, a_r) = a_{\underline{\lambda}} \cdot a_{{\underline{\mu}}^{\star}} \cdot  a_r^{r(g+u) +k},$$ with $u$ the cardinality of the multi-index $\underline{\mu}.$ 
Similarly \eqref{deg2} follows from \eqref{deg4} taking $$P = a_{\underline{\lambda_1}} \cdot a_{\underline {\mu_1}^{\star}} \cdot a_r^{r (g_1+u_1) +k}, \, \, \, Q = a_{\underline {\lambda_2}} a_{\underline{\mu_2}^{\star}} \cdot a_r^{r (g_2 + u_2)},$$ with $u_1, u_2$ being the number of components of $\underline \mu_1, \underline \mu_2$. 

\section*{Appendix: Parabolic Verlinde numbers in arbitrary degree}

In this appendix we check equation \eqref{rgeneral}. Consider a multipartition$$\underline \lambda=(\lambda_p)_{p\in I}$$ labeled by $u$ points $p\in I$ on the curve $C$, such that \begin{equation}\label{nm}kd+|\underline\lambda^{\star}|\equiv 0\mod r.\end{equation} We seek to show that \begin{equation}\label{fgl2}\int_{[Q_C (\mathbb G(r, k+r), d)]^{vir}} a_{\underline \lambda}\cdot a_r^{t}=\chi (U(r, d, \underline \lambda^{\star}), \mathcal L_{\underline \lambda^{\star}}\otimes {\det}^{\star}\Theta_1),\end{equation} where $t$ is chosen so that $$rt+|\underline\lambda|=(r+k)d-rk(g-1).$$ On both sides of \eqref{fgl2} we may consider the degree $\mod r$, hence we may assume $0\leq d<r.$  

The left hand side of \eqref{fgl2} is immediately computed by the Vafa-Intriligator formula \eqref{vaint}: $$(-1)^{d(r-1)}(r+k)^{r\gbar}\sum s_{\underline\lambda}\left(\exp 2 \pi i \frac{ \vec n}{r+k} \right) \exp \left(2\pi i \frac{t-r(g-1)}{r+k}\sum_{i} n_i\right) $$ \begin{equation}\label{vif}\times \prod_{i<j} \left(2\sin \frac{\pi(n_i-n_j)}{r+k}\right)^{-2\gbar}.\end{equation} The sum is taken over all integer vectors $\vec n=(n_1, \ldots, n_r)$ where $$0\leq n_r< \ldots < n_1<r+k.$$ 

We calculate now the right hand side of \eqref{fgl2}. Using the degree $r^{2g}$ tensor product map $$\tau:SU(r, d, \underline \lambda^{\star})\times \text{Jac} \to U(r, d, \underline\lambda^{\star}),$$ we observe the splitting $$\tau^{\star} \mathcal L_{\underline \lambda^{\star}}=\mathcal L_{\underline \lambda^{\star}}\boxtimes \Theta_1^{rk}.$$ Therefore, 
$$\chi(U(r, d, \underline \lambda^{\star}), \mathcal L_{\lambda^{\star}}\otimes {\det}^{\star} \Theta_1)=\left(\frac{r+k}{r}\right)^g\chi (SU\left(r, d, \underline \lambda^{\star}), \mathcal L_{\underline \lambda^{\star}}\right).$$ Write ${\bf SU}(r, d, \underline \lambda^{\star})$ for the moduli stack of parabolic bundles $E$ with determinant isomorphic to $\mathcal O_C(dx)$, for a fixed $x\in C\setminus I$. According to \cite {beauvillelaszlo}, the Euler characteristic can be lifted to the stack, $$\chi(SU\left(r, d, \underline \lambda^{\star}), \mathcal L_{\underline \lambda^{\star}}\right)=\chi({\bf SU}\left(r, d, \underline \lambda^{\star}), q^{\star}\mathcal L_{\underline \lambda^{\star}}\right),$$ via the forgetful map $$q:{\bf SU} (r, d, \underline \lambda^{\star})^{ss}\to SU(r, d, \underline \lambda^{\star}).$$ Here ${\bf SU} (r, d, \underline \lambda^{\star})^{ss} \subset {\bf SU} (r, d, \underline \lambda^{\star})$ denotes the open substack of semistable parabolic bundles.

We write $SL_r(U)$ for $SL_r$ matrices whose entries are regular functions on the open set $U\subset C$. The quotient presentation $${\bf SU}(r, d, \underline \lambda^{\star})=SL_r(C-x)\setminus \mathbf A_{d}^{par}$$ was observed in \cite {pauly}. Here, $\mathbf A_{d}^{par}$ is the ``degree $d$ parabolic affine Grassmannian," which splits as a product $$\mathbf A_d^{par}=\mathbf A_{d}\times \text{Fl}_{\underline \lambda^{\star}}.$$ In turn, $\mathbf A_{d}$ parametrizes pairs $(E, \rho)$ of bundles $E$ and trivializations $\rho$ of $E$ on $C\setminus \{x\}$ which extend to isomorphisms $$\det \rho:\det E\to \mathcal O_C(dx)$$ on the curve $C$. Moreover, $\text{Fl}_{\underline \lambda^{\star}}$ is the product of flag varieties of type $\underline \lambda^{\star}$. The group $SL_r(C-x)$ naturally acts on $\text{Fl}_{\lambda^{\star}}$ via the morphism $$ev_p: SL_r(C-x)\to SL_r, \text { for all } p\in I.$$ It was shown in Remark 3.6 of \cite {beauvillelaszlo}, that $$\mathbf A_{d}=\gamma^{-1} \mathbf A \gamma$$ for some element $\gamma$ in $GL_{r}(\mathbb C((z))$ of order $d$. Here, we write $\mathbf A=\mathbf A_0$ for the degree $0$ affine Grassmannian.

The pullback of the line bundle $\mathcal L_{\underline \lambda^{\star}}$ to $\mathbf A^{par}$ splits as a product $$ \mathcal L_\chi^{k}\otimes_{p\in I} \mathcal N_{\lambda_p^{\star}}$$ The line bundle $\mathcal L_{\chi}^{k}$ over $\mathbf A$ is explicitly described in Lemma 9.2 of \cite {beauvillelaszlo}, and $\mathcal N_{\lambda^{\star}_p}$ are the Borel-Weil bundles on $\text{Fl}_{\underline \lambda^{\star}}$. The space of sections of the product of Borel-Weil bundles is isomorphic to the product $$\mathbf W_{\underline\lambda^{\star}}=\otimes_{p\in I} \mathbf W_{\lambda_p^{\star}}$$ of representations of $\mathfrak{sl}_{r}$ corresponding to the highest weights $\underline \lambda^{\star} = (\lambda_p^{\star})_{p \in I}$. We let ${\mathbf V}_{\lambda}$ denote the representation of $\widehat {\mathfrak {sl}}_{r}$ of level $k$ and highest weight $\lambda.$ The sections of the line bundle $\mathcal L_{\chi}^k$ are shown in Theorem $9.3$ \cite{beauvillelaszlo} to correspond to $\mathbf V_{k\omega}^{\vee}$, where $\omega$ is the fundamental weight $$\omega=e_1+\ldots+e_{r-d},$$ expressed in terms of the standard coordinates.

Write $$J=I\cup \{x\},\, \, \,   \underline\mu=\underline \lambda \cup \{k\omega\}.$$ We regard $\underline \mu$ as a multipartition  indexed by $J$. The $SL_r(C-x)$ invariants in the product $$\mathbf V_{k\omega, k}^{\vee}\otimes {\mathbf W}_{\underline \lambda^{\star}}$$ are isomorphic  to the $SL_r(C-J)$ invariants in ${\mathbf V}_{\underline \mu}^{\vee}$ cf. Proposition $2.3$ \cite{beauville0}. The dimension of the space of invariants is calculated in Corollary $9.8$ of \cite {beauville0}: $$(r+k)^{r\gbar}\left(\frac{r}{r+k}\right)^{g-1}\sum_{\vec\nu}\text{Trace}_{k\omega}\left(\exp 2\pi i
\frac{\vec\nu}{r+k}\right)
\text{Trace}_{\underline \lambda}\left(\exp 2\pi i
\frac{\vec\nu}{r+k}\right)$$ \begin{equation}\label{fgl3}\times\prod_{i<j}\left(2\sin \pi\frac{\nu_i-\nu_j}{r+k}\right)^{-2\gbar}.\end{equation} Here, $\text{Trace}_{\lambda}(g)$ denotes the trace of $g$ in the $SU(r)$ representation of highest weight $\lambda$. In the summation, $\vec\nu$ is an element of the weight lattice of $\mathfrak{sl}_{r}$ whose standard coordinates satisfy the property
$$\nu_1>\ldots> \nu_{r},\, \nu_i-\nu_j \in \mathbb Z,\,\nu_1-\nu_{r}<r+k, \text { and } \nu_1+\ldots+\nu_{r}=0.$$ The first trace was calculated in Proposition $5$ of \cite{sd}: $$(-1)^{d(r-1)}\exp(2\pi i(\nu_1+\ldots+\nu_{r-d})).$$ In fact, we will follow the argument in Proposition $5$ of \cite {sd} from now on, introducing new integer coordinates $$n_i=\nu_i-\nu_{r}.$$ We have $$\nu_{r}= - \frac{1}{r} \sum n_i,$$ hence the trace in the $k\omega$ highest weight representation rewrites $$(-1)^{d(r-1)}\exp\left(2\pi i \frac{d}{r} \sum n_i \right).$$ The traces with respect to $\underline \lambda$ in \eqref {fgl3} equal $$s_{\underline \lambda}\left(\exp 2\pi i \frac{\vec \nu}{r+k}\right)=s_{\underline \lambda}\left(\exp 2\pi i \frac{\vec n}{r+k}\right)\cdot \exp \left (2\pi i \frac{\nu_{r} }{r+k}|\underline \lambda|\right),$$ by the homogeneity of the Schur polynomials. 

Therefore, from \eqref{fgl3} we conclude that, even in the absence of the selection rule \eqref{nm}, the Verlinde Euler characteristic on the moduli stack $$\chi({\bf SU}(r, d, \underline \lambda^{\star}), \mathcal L_{\underline \lambda^{\star}})$$ equals $$(-1)^{d(r-1)}(r+k)^{r \gbar}\left(\frac{r}{r+k}\right)^{g} \sum_{0\leq n_r<\ldots < n_1<r+k}\exp\left(2\pi i\left(\frac{d}{r}-\frac{|\underline\lambda|}{r(r+k)}\right)\sum n_i\right) $$\begin{equation}\label{dimm}\times s_{\underline \lambda} \left(\exp 2\pi i \frac{\vec n}{r+k}\right) \cdot \prod_{i<j}\left(2\sin\pi \frac{n_i-n_j}{r+k}\right)^{-2\gbar}.\end{equation} Furthermore, when \eqref{nm} is satisfied, we have $$\frac{d}{r}-\frac{|\underline\lambda|}{r(r+k)}=\frac{t-r(g-1)}{r+k}+(g-1).$$ Equation \eqref{dimm} then matches the Vafa-Intriligator expression \eqref{vif}, as claimed.

\end{document}